\documentstyle[12pt]{article}

\newtheorem{thm}{Theorem}
\newtheorem{lem}[thm]{Lemma}
\newtheorem{cor}[thm]{Corollary}

\newcommand{\bx}{\rule{.6em}{.6em}}
\newenvironment{pf}{\noindent {\em Proof:}}{\hfill$\bx$\\}

\newcommand{\cz}{c_0}
\newcommand{\kw}{K_{w*}}
\newcommand{\czz}{c_{00}}
\newcommand{\epten}{\tilde{\otimes}_\epsilon}
\newcommand{\ep}{\epsilon}
\newcommand{\la}{\langle}
\newcommand{\ra}{\rangle}
\newcommand{\noteq}{\neq}
\newcommand{\str}{strongly }
\newcommand{\inc}{incomparable }
\newcommand{\incc}{incomparable}
\newcommand{\D}{\Delta}
\newcommand{\N}{\mbox{\hskip.1em N \hskip -1.25em \relax I \hskip
.1em}}
\newcommand{\R}{\mbox{\hskip.1em R \hskip -1.25em \relax I \hskip
.1em}}
\newcommand{\levn}{\{0,1\}^n}
\newcommand{\levjn}{\{0,1\}^{j_n}}
\newcommand{\levjnm}{\{0,1\}^{j_{n_m}}}

\newcommand{\spa}{\hspace{1em}}

\newcommand{\Pinf}{{\cal P}_\infty}
\newcommand{\tdot}{|\la T_{n_m}\gamma,\xi\ra|}
\newcommand{\tdott}{|\la T_{n_m}\tilde{\gamma},\xi\ra|}
\newcommand{\tdotnot}{|\la T_{n_{m_0}}\gamma,\xi\ra|}
\newcommand{\tdotnt}{|\la T_n\tilde{\gamma},\xi\ra|}
\newcommand{\tdotki}{|\la T_{n_m}\gamma_{k,i},\xi_{k,i}\ra|}
\newcommand{\tdotdzp}{\la T_{n_m}\delta_{k,i},\zeta_{k',i'}\ra}
\newcommand{\tdotdz}{\la T_{n_m}\delta_{k,i},\zeta_{k',i'}\ra}

\begin{document}

\begin{center}
  {\Large\bf Some stability properties\\ of $c_0$-saturated
           spaces}\vspace{3mm}\\
  {\large\sc Denny H. Leung}
\end{center}

\vspace{1mm}

\begin{abstract}
A Banach space is $c_0$-saturated if all of its closed infinite
dimensional subspaces contain an isomorph of $c_0$.  In this article,
we study the stability of this property under the formation of
direct sums and tensor products.  Some of the results are:
(1) a slightly more general version of the fact that $c_0$-sums of
$c_0$-saturated spaces are $c_0$-saturated; (2) $C(K,E)$ is
$c_0$-saturated if both $C(K)$ and $E$ are; (3) the tensor product
$JH\tilde{\otimes}_\epsilon JH$ is $c_0$-saturated, where $JH$ is the
James Hagler space.
\end{abstract}

\begin{figure}[b]
  \rule{3in}{.005in}
\end{figure}

Let $E$ be a Banach space.  Following Rosenthal \cite{R},
we say that a Banach
space $F$ is  $E$-saturated if every infinite dimensional closed
subspace of $F$ contains an isomorphic copy of $E$.  In this article,
we will be concerned with the stability properties of $\cz$-saturated
spaces under the formation of direct sums and tensor products. In \S
1, we prove a result which implies that $\cz$-sums of $\cz$-saturated
spaces are $\cz$-saturated.  In \S 2, it is shown that the tensor
product $E \epten F$ is $\cz$-saturated if $E$ is isomorphically
polyhedral (see \S 2 for the definition) and $F$ is $\cz$-saturated.
As a corollary, we obtain that $C(K,E)$ is $\cz$-saturated if and only
if both $C(K)$ and $E$ are.  Finally, in \S 3, we show that $JH \epten
JH$ is $\cz$-saturated, where $JH$ denotes the James Hagler space
\cite{H}.

Standard Banach space terminology, as may be found in \cite{LT}, is
employed. The (closed) unit ball of a Banach space $E$ is denoted by
$U_E$.  The space $\czz$ consists of all finitely non-zero real
sequences.  If $(x_n)$ and $(y_n)$ are sequences residing in possibly
different Banach spaces, we say that $(x_n)$ {\em dominates}\/ $(y_n)$
if there is a constant $K < \infty$ such that $\|\sum a_ny_n\| \leq
K\|\sum a_nx_n\|$ for all $(a_n) \in \czz$.  Two sequences which
dominate each other are said to be {\em equivalent}.  A sequence
$(x_n)$ in a Banach space is {\em semi-normalized}\/ if $0 <
\inf\|x_n\| \leq \sup\|x_n\| < \infty$.
If $A$ is an arbitrary set, $|A|$ denotes the cardinality of
$A$. For an infinite set $A$, $\Pinf(A)$ is the set of all infinite
subsets of $A$.

\section{Direct sums of $\cz$-saturated spaces}

In \cite{R}, it is stated without proof that $\cz$-sums of
$\cz$-saturated spaces are $\cz$-saturated.  In this section, we prove
a result which includes this as a special case.  Let $(E_n)$ be a
sequence of Banach spaces, and let $F$ be a Banach space with a basis
$(e_n)$. The $F$-sum of the spaces $E_n$ is the Banach space $(\oplus
E_n)_F$ of all
sequences $(x(n))$ such that $x(n) \in E_n$ for all $n$, and
$\sum\|x(n)\|e_n$ converges in $F$, endowed with the norm
\[ \|(x(n))\| = \|\sum\|x(n)\|e_n\|. \]
For convenience, we will say that a Banach space is $p$-saturated if
it is $\ell^p$-saturated ($1 \leq p < \infty$) or $\cz$-saturated ($p
= \infty$).

\begin{lem}\label{EF}
Let $E, F$ be $p$-saturated Banach spaces for some $1 \leq p \leq
\infty$, then $E\oplus F$ is $p$-saturated.
\end{lem}

\begin{pf}
It suffices to show that every normalized basic sequence in
$E\oplus F$ has a block basis equivalent to the $\ell^p$-basis
($c_0$-basis if $p = \infty$).
Let $(x_n\oplus y_n)$ be a normalized basic sequence in
$E\oplus F$.  If $p = 1$, and  $(x_n\oplus y_n)$ has a
$\ell^1$-subsequence,
then we are done. Otherwise, by Rosenthal's Theorem \cite{Ro}, we
may assume that
$(x_n\oplus y_n)$ is weakly Cauchy. If $p \neq 1$, using again
Rosenthal's Theorem, we may assume that both $(x_n)$ and
$(y_n)$ are weakly Cauchy.  In both cases, by replacing the sequence
$(x_n\oplus y_n)$
with $(x_{2n-1}-x_{2n}\oplus y_{2n-1}-y_{2n})$ if necessary, we may
even assume that both $(x_n)$ and $(y_n)$ are weakly null.  If
$\|y_n\| \to 0$, then a subsequence of  $(x_n\oplus y_n)$ is
equivalent
to a subsequence of  $(x_n\oplus 0)$.  But then the latter is a
basic sequence in $E$, and hence has a block basis equivalent to the
$\ell^p$-basis.  Therefore,  $(x_n\oplus y_n)$ has a $\ell^p$-block
basis as well.  A similar argument holds if $\|x_n\| \to 0$.
Otherwise, we may take both $(x_n)$ and $(y_n)$ to be semi-normalized
weakly null sequences.  By using a subsequence, it may be assumed that
both are basic sequences.  Then $(x_n)$ has a $\ell^p$-block basis
$(u_k) = (\sum^{n_{k+1}}_{j=n_k+1}a_jx_j)$. Let $(v_k) =
(\sum^{n_{k+1}}_{j=n_k+1}a_jy_j)$ be the corresponding block basis of
$(y_n)$.  If $\|v_k\| \to 0$, then we apply the same argument as
above.  If $\|v_k\| \to \infty$, then
\[ \left(\frac{u_k}{\|v_k\|}\oplus\frac{v_k}{\|v_k\|}\right) \]
is a semi-normalized block basis of $(x_n\oplus y_n)$.  Since
$\|u_k\|/\|v_k\| \to 0$, we may apply the argument above yet again to
conclude the proof.  Finally, then, we may assume that $(v_k)$ is
semi-normalized. Then $(v_k)$ has a $\ell^p$-block basis $(t_k)$.  Let
$(s_k)$ be the corresponding block basis of $(u_k)$ formed by using
the same coefficients.  Arguing as before, we may assume that $(s_k)$
is semi-normalized.  But since $(u_k)$ is a $\ell^p$-sequence, so is
$(s_k)$.
Therefore, $(s_k\oplus t_k)$ is a $\ell^p$-block basis of $(u_k\oplus
v_k)$, and hence of $(x_n\oplus y_n)$.
\end{pf}

\begin{thm}
Let $(E_n)$ be a sequence of $p$-saturated Banach spaces, and let $F$
be a
$p$-saturated Banach space with a basis.  Then $E = (\oplus E_n)_F$ is
$p$-saturated.
\end{thm}

\begin{pf}
For each $x \in E$, write $x = (x(n))$, where $x(n) \in E_n$ for all
$n$.  Let $(x_k)$ be a normalized basic sequence in $E$. For each $m
\in \N$, let $P_m$ be
the projection on $E$ defined by $P_mx = y$, where $y(n) = x(n)$ if $n
\leq m$, and $y(n) = 0$ otherwise.  If for some subsequence $(z_j)$ of
$(x_k)$, and some $m \in \N$,  $(P_mz_j)_j$ dominates
$(z_j)$, then $(z_j)$ is equivalent to $(P_mz_j)_j$.  But then
the latter is a basic sequence in $E_1\oplus\ldots\oplus E_m$, which
is $p$-saturated by Lemma \ref{EF} and induction.  Hence $(z_j)$, and
thus
$(x_k)$, has a $\ell^p$-block basis, and we are done.  Otherwise, for
all $m \in \N$, and every subsequence $(z_j)$ of $(x_k)$,
\begin{equation}\label{nodom}
\inf\{\|\sum_j a_jP_mz_j\| : (a_j) \in c_{00},\ \|\sum a_jz_j\| =
1\} = 0 .
\end{equation}
Let $m_0 = 0$, and $y_1 = x_1$. Choose $m_1$ such that
$\|(1-P_{m_1})y_1\| \leq 1$.   By
(\ref{nodom}), there exists $k_2 \geq 2$, and $y_2 \in
\mbox{span}\{x_k: 2 \leq k \leq k_2\}$, $\|y_2\| = 1$, such that
$\|P_{m_1}y_2\| \leq
1/4$.  Then choose $m_2 > m_1$ so that $\|(1-P_{m_2})y_2\| \leq 1/4$.
Continuing inductively, we obtain a normalized block basis $(y_j)$ of
$(x_k)$,
and $(m_j)^\infty_{j=0}$ such that
\[ \|y_j - (P_{m_j}-P_{m_{j-1}})y_j\| \leq 1/j \]
for all $j \geq 1$  ($P_0 = 0$).  Let $v_j =
(P_{m_j}-P_{m_{j-1}})y_j$.  Then $(y_j)$ has a subsequence equivalent
to a subsequence of $(v_j)$.  But, writing the basis of $F$ as
$(e_n)$, it is clear that $(v_j)$ is equivalent to the sequence
$(\sum^{m_j}_{n=m_{j-1}+1}\|v_j(n)\|e_n)$ in $F$.  Since $F$ is
$p$-saturated,  we conclude that any subsequence of $(v_j)$ has a
$\ell^p$-block basis.
Thus, the same can be said of $(x_k)$, and the proof is complete.
\end{pf}

\section{Tensor products of $\cz$-saturated spaces}

For Banach spaces $E$ and $F$, let $\kw(E',F)$ denote the space of all
compact weak*-weakly continuous operators from $E'$ into $F$, endowed
with the operator norm.  The
$\ep$-tensor product $E \epten F$ is the closure in $\kw(E',F)$ of the
finite rank operators that belong to $\kw(E',F)$.  These spaces are
equal if either $E$ or $F$ has the approximation property \cite{Ru}.
In this
section and the next, we investigate special cases of the following:
\\

\noindent{\em Problem}: Is $\kw(E',F)$ (or $E\epten F$)
$\cz$-saturated if both $E$ and $F$ are? \\

A Banach space is {\em polyhedral}\/ if the unit ball of every finite
dimensional subspace is a polyhedron.  It is {\em isomorphically
polyhedral}\/ if it is isomorphic to a polyhedral Banach space.  Our
interest in isomorphically polyhedral spaces arises from the following
result of Fonf \cite{F1}.

\begin{thm}
{\em (Fonf)} An isomorphically polyhedral Banach space is
$\cz$-saturated.
\end{thm}

A subset $W$ of the dual of a Banach space $E$ is said to be {\em
isomorphically precisely norming}\/ (i.p.n.)\ if $W$ is bounded and\\
(a) there exists $K < \infty$ such that $\|x\| \leq K\sup_{w\in
W}|w(x)|$ for all $x \in E$, \\
(b) the supremum $\sup_{w\in W}|w(x)|$ is attained at some $w_0 \in W$
for all $x \in E$. \\
This terminology was introduced by Rosenthal \cite{R1,R} to provide a
succint formulation of the following result of Fonf \cite{F2}.

\begin{thm}
{\em (Fonf)} A separable Banach space $E$ is isomorphically
polyhedral if and only if $E'$ contains a countable i.p.n.\ subset.
\end{thm}

In this section, we consider the space $\kw(E',F)$ when one of the
spaces $E$ or $F$ is isomorphically polyhedral, and the other is
$\cz$-saturated.  Note the symmetry in the situation as $\kw(E',F)$ is
isometric to $\kw(F',E)$ via the mapping $T \mapsto T'$. For Lemma
\ref{pn}, note that if $x' \in E'$ and $y' \in F'$, the pair $(x',y')$
defines a functional on $\kw(E',F)$ by $T \mapsto \la Tx',y'\ra$.

\begin{lem}\label{pn}
Let $E, F$ be Banach spaces, and let $W$ and $V$ be i.p.n.\
subsets of $E'$ and $F'$ respectively.  Then
$W\times V$ is an i.p.n.\ subset of
$(\kw(E',F))'$.
\end{lem}

\begin{pf}
It is clear that if both $W$ and $V$ satisy (a) of the definition of
an i.p.n.\ set with constant $K$, then $W\times V$ also satisfies it
with constant $K^2$.  Now assume that both $W$ and $V$ satisfy part
(b) of the definition.  It is easy to see that $(x',y') \mapsto \la
Tx',y'\ra$ is a continuous function on
$\overline{W}^{w^*}\times\overline{V}^{w^*}$, where both
$\overline{W}^{w^*}$ and $\overline{V}^{w^*}$ are given their
respective weak* topologies.  Since
$\overline{W}^{w^*}\times\overline{V}^{w^*}$ is compact, there exists
$(w_0,v_0) \in \overline{W}^{w^*}\times\overline{V}^{w^*}$ such that
\[ \sup_{(x',y')\in W\times V}|\la Tx',y'\ra| =
\sup_{(x',y')\in\overline{W}^{w^*}\times\overline{V}^{w^*}}|\la
Tx',y'\ra| = |\la Tw_0,v_0\ra| . \]
Now there exists $v\in V$  such that
\begin{eqnarray*}
|\la Tw_0,v\ra| & = & \sup_{y'\in V}|\la Tw_o,y'\ra| \\
 & = & \sup_{y'\in\overline{V}^{w^*}}|\la Tw_o,y'\ra| \\
 & \geq & |\la Tw_0,v_0\ra|.
\end{eqnarray*}
Similarly, there exists $w \in W$ such that
\begin{eqnarray*}
|\la Tw,v\ra| & = & |\la w,T'v\ra| \\
& = & \sup_{x'\in W}|\la x',T'v\ra| \\
& = & \sup_{x'\in\overline{W}^{w^*}}|\la x',T'v\ra| \\
& \geq & |\la Tw_0,v\ra|.
\end{eqnarray*}
Combining the above, we see that $|\la Tw,v\ra| \geq
\sup_{(x',y')\in W\times V}|\la Tx',y'\ra|$.  Since the reverse
inequality is obvious, the proof is complete.
\end{pf}

\begin{lem}\label{block}
Let $(x_n)$ be a non-convergent sequence in a $\cz$-saturated Banach
space $F$. There exists a normalized block $(u_k) =
(\sum^{n_k}_{n=n_{k-1}+1}b_nx_n)$ of $(x_n)$ which is equivalent to
the $\cz$-basis.
\end{lem}

\begin{pf}
Going to a subsequence, we may assume that $\inf_{m\neq n}\|x_m-x_n\|
> 0$.
By Rosenthal's Theorem \cite{Ro}, we may also assume that $(x_n)$ is
weakly
Cauchy. Let $y_n = (x_{2n-1}-x_{2n})/\|x_{2n-1}-x_{2n}\|$ for all $n$.
Then $(y_n)$ is a weakly null normalized  block of $(x_n)$.  Without
loss of generality, we may assume that $(y_n)$ is a basic sequence.
Since $F$ is $\cz$-saturated, $(y_n)$ has a normalized block basis
$(u_k)$ equivalent to the $\cz$-basis. As $(u_k)$ is also a normalized
block of $(x_n)$, the proof is complete.
\end{pf}

\begin{lem}\label{opbl}
Let $E, F$ be Banach spaces so that $F$ is $\cz$-saturated, and let
$(T_n)$  be a normalized basic sequence in $\kw(E',F)$.
For every $x' \in E'$, there is a normalized block
basis $(S_n)$ of $(T_n)$, and a constant $C$, such that
\[ \|\sum a_nS_nx'\| \leq C\sup_n|a_n| \]
for all $(a_n) \in \czz$.
\end{lem}

\begin{pf}
There is no loss of generality in assuming that $\|x'\| = 1$.\\

\noindent\underline{Case 1}  $(T_nx')$ converges.\\
We may assume without loss of generality that $\|(T_{2n-1}-T_{2n})x'\|
\leq 2^{-n}$ for all $n$.  Let $\epsilon_n = \|T_{2n-1}-T_{2n}\|$.
Since $(T_n)$ is a normalized basic sequence, $\epsilon \equiv
\inf_n\epsilon_n > 0$. Now let $S_n =
\epsilon^{-1}_n(T_{2n-1}-T_{2n})$ for all $n$.  Then $(S_n)$ is a
normalized block basis of $(T_n)$.  Furthermore, for any $(a_n) \in
\czz$,
\begin{eqnarray*}
\|\sum a_nS_nx'\| & \leq &
\epsilon^{-1}\sup_n|a_n|\sum\|(T_{2n-1}-T_{2n})x'\| \\
& \leq & \epsilon^{-1}\sup_n|a_n| .
\end{eqnarray*}
\underline{Case 2} $(T_nx')$ does not converge.  \\
By Lemma \ref{block}, $(T_nx')$ has a normalized block $(u_k) =
(\sum^{n_k}_{n=n_{k-1}+1}b_nT_nx')$ which is equivalent to the
$\cz$-basis.  Let $R_k = (\sum^{n_k}_{n=n_{k-1}+1}b_nT_n)$.  Then
$\|R_k\| \geq \|R_kx'\| = \|u_k\| = 1$.  Let $S_k = R_k/\|R_k\|$ for
all $k$.  Then $(S_k)$ is a normalized block basis of $(T_n)$, and
\[
\|\sum a_kS_kx'\|
 =  \|\sum \frac{a_k}{\|R_k\|}u_k\| \leq C\sup_k|a_k| \]
for all $(a_k) \in \czz$, since $(u_k)$ is equivalent to the
$\cz$-basis and $\|R_k\| \geq 1$.
\end{pf}

\begin{thm}\label{main}
Let $E, F$ be Banach spaces so that $E$ is isomorphically polyhedral
and $F$ is $\cz$-saturated.  Then $\kw(E',F)$ is $\cz$-saturated.
\end{thm}

\begin{pf}
Let $(T_n)$ be a normalized basic sequence in $\kw(E',F)$.  It is
easily seen that $G = [\cup T'_nF']$ is a separable subspace of $E$,
and that the sequences $(T_n)$ and $(T_{n|G})$ are equivalent.  Thus
we may assume that $E$ is separable.  Then $E'$ contains a countable
i.p.n.\ subset $W$.  Write $W = (w_m)$. By
Lemma \ref{opbl}, $(T_n)^\infty_{n=1}$ has a normalized block basis
$(T^{(1)}_n)$ such that
\[ \|\sum a_nT^{(1)}_nw_1\| \leq C_1\sup_n|a_n| \]
for some constant $C_1$ for all $(a_n) \in \czz$.  Inductively, if
$(T^{(m)}_n)$ has been chosen, let $(T^{(m+1)}_n)$ be a normalized
block basis of $(T^{(m)}_n)^\infty_{n=2}$ such that there is a
constant $C_{m+1}$ satisfying
\[ \|\sum a_nT^{(m+1)}_nw_{m+1}\| \leq C_{m+1}\sup_n|a_n| \]
for all $(a_n) \in \czz$.  Now let $S_m = T^{(m)}_1$ for all $m \in
\N$.  Then $(S_m)$ is a normalized block basis of $(T_n)$.  Also, for
all $k$, $(S_m)^\infty_{m=k}$ is a block basis of $(T^{(k)}_n)$.  Fix
$k$, and write $S_m = \sum^{j_m}_{n=j_{m-1}+1}b_nT^{(k)}_n$ for all $m
\geq k$. Since $(S_m)$ is normalized, and $(T^{(k)}_n)$ is normalized
basic, $(b_i)$ is bounded.  Therefore, by the choice of $(T^{(k)}_n)$,
\begin{eqnarray*}
\left\|\sum^\infty_{m=k}a_mS_mw_k\right\| & = &
\left\|\sum^\infty_{m=k}a_m\left(\sum^{j_m}_{n=j_{m-1}+1}
b_nT^{(k)}_nw_k\right)\right\| \\
& \leq & C_k\sup_{m\geq k}\sup_{j_{m-1}<n\leq j_m}|a_mb_n| \\
& \leq & C_k\sup_n|b_n|\sup_{m\geq k}|a_m|
\end{eqnarray*}
for all $(a_m) \in \czz$.  Consequently,
$\sum^\infty_{m=1}a_mS_mw_k$ converges in $F$ for all $k$ and all
$(a_m)
\in \cz$.  Hence $\sum|\la S_mw,y'\ra|
< \infty$ for all $(w,y') \in W\times U_{F'}$.  But by Lemma \ref{pn},
$W\times U_{F'}$ is an i.p.n.\ subset of
$\kw(E',F)$. Applying Elton's extremal criterion (\cite{E}, see also
\cite[Theorem 18]{R}), we see that $[S_m]$ contains a copy of $c_0$.
\end{pf}

Recall that a subset $A$ of topological space $X$ is
{\em dense-in-itself}\/ if every point of $A$ is an accumulation point
of $A$.  $A$ is {\em scattered}\/ if it contains no non-empty
dense-in-itself subset. \\

\begin{cor}
Let $K$ be a compact Hausdorff space, and let $E$ be a Banach space.
Then $C(K,E)$ is $\cz$-saturated if and only if both $C(K)$ and $E$
are $\cz$-saturated.
\end{cor}

\begin{pf}
The ``only if'' part is clear, since both $C(K)$ and $E$ embed in
$C(K,E)$.  Now assume that both $C(K)$ and $E$ are $\cz$-saturated.
To begin with, assume additionally that $C(K)$ is separable.  Then $K$
is metrizable \cite[Proposition II.7.5]{S}. If $K$
is not scattered, by \cite[Theorem 8.5.4]{Sem}, there is a continuous
surjection $\phi$ of $K$ onto $[0,1]$.  Then $f \mapsto f\circ\phi$ is
an isometric embedding of $C[0,1]$ into $C(K)$.  This contradicts the
fact that $C(K)$ is $\cz$-saturated.  Thus $K$ is scattered.  By
\cite[Theorem 8.6.10]{Sem}, $K$ is homeomorphic to a countable compact
ordinal.  In particular, $K$ is countable.  Hence $C(K)'$ contains a
countable i.p.n.\ subset, namely, $\{\delta_k: k \in
K\}$, where $\delta_k$ denotes the Dirac measure concentrated at $k$.
Therefore, $C(K)$ is isomorphically polyhedral, and $C(K,E) =
\kw(C(K)',E)$ is $\cz$-saturated by Theorem \ref{main}.

If $C(K)$ is non-separable, as in the proof of the
theorem, it suffices to show that $\kw(G',E)$ contains a copy of $\cz$
for an arbitrary separable closed subspace $G$ of $C(K)$.  However,
$\kw(G',E)$ is isometric to $\kw(E',G)$, which clearly embeds in
$\kw(E',F)$ for any closed subspace $F$ of $C(K)$ containing $G$.
Take
$F$ to be the closed sublattice generated by $G$ and the constant $1$
function.  By Kakutani's Representation Theorem \cite[Theorem
II.7.4]{S},
$F$ is lattice isometric to some $C(H)$.  Note that $C(H)$ is
separable since $F$ is.  Therefore, $\kw(E',F)$, which is isometric to
$\kw(F',E) = C(H,E)$, is $\cz$-saturated by the above.  Since
$\kw(G',E)$ is isomorphic to a subspace of $\kw(E',F)$, it contains a
copy of $\cz$.
\end{pf}

\section{The space $JH\epten JH$}

In view of Theorem \ref{main} in \S 2, it is interesting to consider
spaces $\kw(E',F)$ where both $E$ and $F$ are $\cz$-saturated, but
neither is isomorphically polyhedral.  In this section, we investigate
one such case. Namely, when $E = F = JH$, the James Hagler space.
In \cite{L}, it was shown that $JH\epten JH = \kw(JH',JH)$
does not contain an isomorph of $\ell^1$.
We show here that, in fact, $JH\epten JH$
is $\cz$-saturated.  The proof uses Elton's extremal criterion for
weak unconditional convergence, and the
``diagonalization technique'' employed by Hagler to show that every
normalized weakly null sequence in $JH$ has a $\cz$-subsequence.

Let us recall the definition of the space JH, as well as fix some
terms and notation.  Let $T = \cup^\infty_{n=0}\{0,1\}^n$ be the dyadic
tree. The elements of $T$ are called {\em nodes}.  If $\phi$ is a node
of the form $(\ep_i)^n_{i=1}$, we say that $\phi$ has {\em length}\/
$n$ and write $|\phi| = n$.  The length of the empty node is defined
to be $0$.  For $\phi, \psi \in T$ with $\phi = (\ep_i)^n_{i=1}$ and
$\psi = (\delta_i)^m_{i=1}$, we say that $\phi \leq \psi$ if $n \leq
m$ and $\ep_i = \delta_i$ for $1 \leq i \leq n$.  The empty node is
$\leq \phi$ for all $\phi \in T$.  We write $\phi < \psi$ if $\phi
\leq \psi$ and $\phi \neq \psi$. Two nodes $\phi$ and $\psi$ are
{\em incomparable}\/ if neither $\phi \leq \psi$ nor $\psi \leq \phi$
hold. If $\phi \leq \psi$, we say that
$\phi$ is an {\em ancestor}\/ of $\psi$, while $\psi$ is a {\em
descendant}\/ of $\phi$. For any $\phi \in T$, let $\Delta_\phi$ be
the set of all descendants of $\phi$.  If $\phi \leq \psi$, let
\[ S(\phi,\psi) = \{\xi: \phi \leq \xi \leq \psi\}. \]
A set of the form $S(\phi,\psi)$ is called a {\em segment}, or more
specifically, a $m$-$n$ {\em segment}\/ provided $|\phi| = m,$ and
$|\psi| =
n$. A {\em branch}\/ is a maximal totally ordered subset of $T$.
The set of all branches is denoted by $\Gamma$.  A branch $\gamma$
(respectively, a segment $S$) is
said to {\em pass through}\/ a node $\phi$ if $\phi \in \gamma$
(respectively, $\phi \in S$). If
$x: T \to \R$\ is
a finitely supported function and $S$ is a segment, we define (with
slight abuse of notation) $Sx = \sum_{\phi\in S}x(\phi)$.  Similarly,
if $\gamma \in \Gamma$, we define $\gamma(x) = \sum_{\phi\in
\gamma}x(\phi)$.  A set of segments $\{S_1,\ldots,S_r\}$ is {\em
admissible}\/ if they are pairwise disjoint, and there are $m, n \in
\N\cup\{0\}$ such that each $S_i$ is a $m$-$n$ segment.  The James
Hagler
space $JH$ is defined as the completion of the set of all finitely
supported functions $x: T \to \R$\ under the norm:
\[ \|x\| = \sup\left\{\sum^r_{i=1}|S_ix| : S_1,\ldots,S_r \mbox{ is an
admissible set of segments}\right\}. \]
Clearly, all $S$ and $\gamma$ extend to norm $1$ functionals on $JH$.
Finally, if $x: T \to \R$\ is finitely supported, and $n \geq 0$, let
$P_nx: T \to \R$\ be defined by
\[ (P_nx)(\phi) = \left\{ \begin{array}{ll}
                          x(\phi) & \mbox{if $|\phi| \geq n$} \\
                          0      & \mbox{otherwise.}
                          \end{array}
                  \right. \]
Obviously, $P_n$ extends uniquely to a norm $1$ projection on $JH$,
which we denote again by $P_n$.

We begin with some lemmas on ``node management''. Let $\phi$ and
$\psi$ be nodes.  Denote by $A(\phi,\psi)$ denote the
unique node of maximal length such that $A(\phi,\psi) \leq \phi$ and
$A(\phi,\psi) \leq \psi$.  A sequence of nodes $(\phi_n)$  is a
{\em strongly incomparable sequence}\/ if\\
(a) $\phi_n$ and $\phi_m$ are incomparable if $n \neq m$,\\
(b) no family of admissible segments passes through more than two of
the $\phi_n$'s.\\
The first lemma is due to Hagler \cite[Lemma 2]{H}.

\begin{lem}
{\em (Hagler)} Let $(\phi_n)$ be a sequence of nodes with strictly
increasing
lengths.  Then there exists $N \in \Pinf{\em(\N)}$ such that either
$(\phi_n)_{n\in N}$ determines a unique branch of $T$, or
$(\phi_n)_{n\in N}$ is a strongly incomparable sequence.
\end{lem}

\begin{lem}\label{ancestor}
Let $(\phi(n))$ be a strongly \inc sequence of nodes such that
$|\phi(n)| < |\phi(n+1)|$ for all $n$.  Then for all $m \geq n \geq
1$,
$\phi(m) \in \D_{A(\phi(n),\phi(n+1))}$.
\end{lem}

\begin{pf}
Otherwise, there exist $m \geq n \geq 1$ such that $\phi(m) \notin
\D_{A(\phi(n),\phi(n+1))}$. In particular, note that $m \geq n + 2$.
Let $\phi_1 = \phi(n)$, and let $\phi_2$ and $\phi_3$ be the ancestors
of $\phi(n+1)$ and $\phi(m)$ respectively of length
$|\phi(n)|$.  Then $\phi_1
\noteq \phi_2$ since $\phi(n)$ and $\phi(n+1)$ are \incc .   Also
$\phi_1 \noteq \phi_3$ and $\phi_2 \noteq \phi_3$ since $\phi_1,
\phi_2 \in \D_{A(\phi(n),\phi(n+1))}$, while $\phi_3 \notin
\D_{A(\phi(n),\phi(n+1))}$. Let $\psi_1, \psi_2$ be nodes of length
$|\phi(m)|$ which are $\geq \phi(n)$ and $\phi(n+1)$ respectively, and
let $\psi_3 = \phi(m)$.  Then $\{S(\phi_i,\psi_i) : i = 1, 2, 3\}$ is
an admissible set of segments.  However, $\phi(n) \in
S(\phi_1,\psi_1), \phi(n+1) \in S(\phi_2,\psi_2)$, and $\phi(m) \in
S(\phi_3,\psi_3)$, violating the strong incomparability of
$(\phi(k))$.
\end{pf}

\begin{lem}\label{inc}
Let $m \in$ {\em \N}, and let $(\phi(i,j))^\infty_{i=1}$ be a \str
\inc
sequence of nodes for $1 \leq j \leq m$.  Assume that $|\phi(i,j)| =
|\phi(i,j')| < |\phi(i+1,j)|$ whenever $i, j , j' \in {\em \N}$, $1
\leq j,
j' \leq m$, and that for each $i$, $\{\phi(i,j) : 1 \leq j \leq m\}$
are
pairwise distinct.  Then there exists $k_0$ such that $\{\phi(i,j) : i
\geq k_0, 1 \leq j \leq m\}$ are pairwise \incc .
\end{lem}

\begin{pf}
Induct on $m$.  If $m = 1$, there is nothing to prove.  Now assume
that the statement is true for $m - 1$ $(m \geq 2)$.  Let $\{\phi(i,j)
: i \geq 1, 1 \leq j \leq m\}$ be as given.  Without loss of
generality, we
may assume that $\{\phi(i,j) : i \geq 1, 1 \leq j \leq m-1 \}$ are
pairwise \incc.  First observe that if $\phi(i_1,j_1) < \phi(i_2,j_2)$
for some $j_1, j_2 \leq m$, then $A(\phi(i_2,j_2),\phi(i_2+1,j_2))
\geq \phi(i_1,j_1)$.  Indeed, since $A(\phi(i_2,j_2),\phi(i_2+1,j_2))$
and $\phi(i_1,j_1)$ share the same descendant $\phi(i_2,j_2)$, they
are comparable. Hence if the claim fails,
$A(\phi(i_2,j_2),\phi(i_2+1,j_2))
< \phi(i_1,j_1)$.  But then the ancestor of $
\phi(i_2+1,j_2)$ of length $|\phi(i_1,j_1)|$, $\phi(i_1,j_1)$, and
$\phi(i_1,j_2)$, are
distinct.  From this
it is easy to construct an admissible set of segments which pass
through all three nodes $\{\phi(i_1,j_2), \phi(i_2,j_2),
\phi(i_2+1,j_2)\}$, in  violation of their strong incomparability.
This proves the claim.  In particular, under the circumstances, Lemma
\ref{ancestor} implies
$\phi(i,j_2) \in \D_{\phi(i_1,j_1)}$ for all $i \geq i_2$.  The
remainder of the proof is divided into two cases.\\

\noindent\underline{Case 1}\spa There exist $j_1 < m$ and $i_1, i_2
\in \N\ $ such that $\phi(i_1,j_1) \leq  \phi(i_2,m)$. \\
\indent Note that since $\{\phi(i_1,j) : j \leq m\}$ are pairwise
distinct, we must have $\phi(i_1,j_1) < \phi(i_2,m)$.  By the
observation above, we obtain that $\phi(i,m) \in \D_{\phi(i_1,j_1)}$
for all $i \geq i_2$.  However, for any $i' \geq i_2, j < m$,
$\phi(i',j)$ is \inc with $\phi(i_1,j_1)$ by induction.  Hence
$\phi(i',j) \notin \D_{\phi(i_1,j_1)}$.  Thus $\phi(i,m)$ and
$\phi(i',j)$ are \inc whenever $i, i' \geq i_2$ and $j < m$.  This is
enough to show that $\{\phi(i,j) : i \geq i_2, 1 \leq j \leq m\}$ are
pairwise \incc.\\

\noindent\underline{Case 2}\spa For all  $j_1 < m$ and $i_1, i_2
\in \N$, $\phi(i_1,j_1) \not\leq  \phi(i_2,m)$. \\
\indent Let $I = \{i :$ there exist $i' \in \N, j < m$ such that
$\phi(i,m) < \phi(i',j)\}$.  Let $i_1$ and $i_2$ be distinct elements
of $I$.  Choose $i'_1, i'_2 \in \N$, $j_1, j_2 < m$ such that
$\phi(i_k,m) < \phi(i'_k,j_k)$, $k = 1, 2$.  By the observation above,
$\phi(i,j_k) \in \D_{\phi(i_k,m)}$ for all $i \geq i'_k$, $k = 1, 2$.
Now $\phi(i_1,m)$ and $\phi(i_2,m)$ are \inc by assumption.
Therefore, $\D_{\phi(i_1,m)} \cap \D_{\phi(i_2,m)} = \emptyset$.
Combined with the above, we see that $j_1 \neq j_2$.  It follows that
$|I| \leq m-1$.  Now choose $k_0$ such that $i < k_0$ for all $i \in
I$.  By the case assumption, $\phi(i,m)$ is \inc with $\phi(i',j)$
whenever $i \geq k_0$ and $j < m$.  This is enough to show that
$\{\phi(i,j) : i \geq k_0, 1 \leq j \leq m\}$ are
pairwise \incc.
\end{pf}

\begin{lem}\label{seven}
Let $(x_n)$ be a bounded weakly null sequence in $JH$ so that
there
is a  sequence $0 = j_1 < j_2 < j_3 < \cdots$ with
$x_n \in (P_{j_n}-P_{j_{n+1}})JH$ for all $n$.  Assume that
$\sup_{\xi\in\Gamma}|\la x_n,\xi\ra| \leq \epsilon$ for all $n$.
Then there is a subsequence $(x_{n_k})$ such that
\[ \sup_{\xi\in\Gamma}\sum^\infty_{k=1}|\la x_{n_k},\xi\ra| \leq
7\epsilon.  \]
\end{lem}

\begin{pf}
For $m, n \in$ \N, let
\begin{eqnarray*}
F(n,m) & = & \{\phi \in \{0,1\}^{j_n} : \mbox{there exists at least
one branch } \xi \mbox{ through } \phi \mbox{ with }
\\ & & |\la x_n,\xi\ra| >\ep/2^m,
\mbox{ and for all branches $\xi$ through } \phi,\ |\la
x_n,\xi\ra| \leq \ep/2^{m-1} \}.
\end{eqnarray*}
Since $(x_n)$ is bounded, $\sup_n|F(n,m)| < \infty$ for all $m$.
Let $N_1 \in \Pinf(\N)$ be such that $(|F(n,1)|)_{n\in N_1}$ is
constant,
say, $b_1$.  Write $F(n,1) = \{\phi(n,1,i) : 1 \leq i \leq b_1\}$ for
all $n \in N_1$.  Choose $N'_1 \in \Pinf(N_1)$ so that for $1 \leq i
\leq
b_1$, $(\phi(n,1,i))_{n\in N'_1}$ is either a \str \inc sequence or
determines a
branch.  Let
\[ I_1 = \{i \leq b_1 : (\phi(n,1,i))_{n\in N'_1} \mbox{
determines a branch}\}, \]
and let $\Gamma_1$ be the set of branches
determined by some $(\phi(n,1,i))_{n\in N'_1}$ for some $i \in I_1$.
Let $L_1 = \{1,\ldots,b_1\}\backslash I_1.$  By Lemma \ref{inc}, there
exists $N''_1 \in \Pinf(N'_1)$ such that $\{\phi(n,1,i) : n \in N''_1,
i \in L_1\}$ are pairwise incomparable.  Finally, since $\Gamma_1$ is
finite, there exists $n_1 \in N''_1$ such that $|\la
x_{n_1},\gamma\ra| \leq \ep/2$ for all $\gamma \in \Gamma_1$.
Continue inductively to obtain $n_1 < n_2 < \cdots$, numbers $b_1,
b_2, \ldots$, and sets $I_1, I_2,
\ldots$, $L_1, L_2, \ldots$, and $\Gamma_1, \Gamma_2, \ldots$ so that
\begin{enumerate}
\item for all $k \geq m$, $|F(n_k,m)| = b_m$,
\item for all $k \geq m$, $F(n_k,m) = \{\phi(n_k,m,i) : i \leq b_m\}$,
where $(\phi(n_k,m,i))^\infty_{k=m}$ determines a branch if $i \in
I_m$, and is a \str
\inc sequence if $i \in L_m = \{1,\ldots,b_m\}\backslash I_m$,
\item for all $m$, $\Gamma_m$ is the set of branches determined by
$(\phi(n_k,m,i))^\infty_{k=m}$  for some $i \in I_m$,
\item for all $m$, $\{\phi(n_k,m,i) : m \leq k, i \in L_m\}$ are
pairwise \incc,
\item for all $k$, $|\la x_{n_k},\gamma\ra| \leq \ep/2^k$ for every
$\gamma \in \Gamma_1\cup\ldots\cup\Gamma_k$.
\end{enumerate}
For all $k$, let $G(k) =
\{0,1\}^{j_{n_k}}\backslash\cup^k_{m=1}F(n_k,m)$.  Then
$\{F(n_k,1),\ldots,F(n_k,k),G(k)\}$ is a partition of
$\{0,1\}^{j_{n_k}}$ for each $k$.  Fix $\xi \in \Gamma$.  Say $\xi =
(\phi_n)^\infty_{n=0}$, where $|\phi_n| = n$ for all $n \geq 0$.  Now
let $J_0 = \{k : \phi_{n_k} \in G(k)\}$, and let $J_m = \{k \geq m :
\phi_{n_k} \in F(n_k,m)\}$ for all $m \geq 1$.  Then
$(J_m)^\infty_{m=0}$ is a partition of \N.  For all $k$ such that $k
\in J_m$ for some $m \geq 1$, choose $i_k$ such that $\phi_{n_k} =
\phi(n_k,m,i_k)$.  For all $m \geq 1$, let $J_{m,1} = \{k \in J_m :
i_k \in I_m\}$, and let $J_{m,2} = \{k \in J_m : i_k \in L_m\}$.  Fix
$m \geq 1$ such that $J_{m,1} \neq \emptyset$.\\

\noindent\underline{Case 1}\spa $J_{m,1}$ is infinite.\\
\indent In this case, there exists $\gamma \in \Gamma_m$ such that
$\gamma = \xi$.  Therefore,
\begin{eqnarray*}
\sum_{k\in J_{m,1}}|\la x_{n_k},\xi\ra| & \leq & \sum^\infty_{k=m}|\la
x_{n_k},\xi\ra| = \sum^\infty_{k=m}|\la x_{n_k},\gamma\ra| \\
& \leq & \sum^\infty_{k=m}\frac{\ep}{2^k} = \frac{\ep}{2^{m-1}}.
\end{eqnarray*}

\noindent\underline{Case 2}\spa $J_{m,1}$ is finite. \\
\indent Let $k_0 = \max J_{m,1}$.  There exists $\gamma \in \Gamma_m$
such that $\phi_{n_{k_0}} \in \gamma$.  Then, since $\phi_{n_{k_0}}
\in
F(n_{k_0},m)$,
\begin{eqnarray*}
\sum_{k\in J_{m,1}}|\la x_{n_k},\xi\ra| & \leq & \sum^{k_0}_{k=m}|\la
x_{n_k},\xi\ra| = \sum^{k_0-1}_{k=m}|\la x_{n_k},\gamma\ra| +
|\la x_{n_{k_0}},\xi\ra| \\
 & \leq &  \sum^\infty_{k=m}\frac{\ep}{2^k} + \frac{\ep}{2^{m-1}} =
\frac{\ep}{2^{m-2}}.
\end{eqnarray*}
Therefore, in either case, we have
\begin{equation}\label{jmone}
\sum_{k\in J_{m,1}}|\la x_{n_k},\xi\ra| \leq
\frac{\ep}{2^{m-2}}.
\end{equation}
Now suppose for some $m \geq 1$, there are distinct $k_1, k_2 \in
J_{m,2}$.  Then $k_1, k_2 \geq m$, and $\phi_{n_{k_l}} =
\phi(n_{k_l},m,i_{k_l})$, for some $i_{k_l} \in L_m$, $l = 1, 2$. But
then by choice, $\phi_{n_{k_1}}$ and $\phi_{n_{k_2}}$ must be \incc.
This is a contradiction since they both belong to the branch $\xi$.
Thus, for all $m \geq 1$, $|J_{m,2}| \leq 1$. Now $k \in J_{m,2}$
implies $\phi_{n_k} \in F(n_k,m)$, and hence $|\la x_{n_k},\xi\ra|
\leq \ep/2^{m-1}$. Consequently, for all $m \geq 1$
\begin{equation}
\sum_{k\in J_{m,2}}|\la x_{n_k},\xi\ra| \leq \frac{\ep}{2^{m-1}}.
\end{equation}
Finally,
\begin{equation}\label{jz}
\sum_{k\in J_0}|\la x_{n_k},\xi\ra| \leq
\sum^\infty_{k=1}\frac{\ep}{2^k} = \ep.
\end{equation}
Combining  inequalities (\ref{jmone})--(\ref{jz}), we obtain
\begin{eqnarray*}
\sum_k|\la x_{n_k},\xi\ra| & \leq &
\sum^\infty_{m=1}\sum_{k\in J_m}|\la x_{n_k},\xi\ra| +
\sum_{k\in J_0}|\la x_{n_k},\xi\ra| \\
& \leq & \sum^\infty_{m=1}\frac{3\ep}{2^{m-1}} + \ep = 7\ep.
\end{eqnarray*}
\end{pf}

For $n \geq 0$, call a subset $D$ of $\levn \times \levn$
{\em diagonal}\/ if whenever
 $(\phi_1, \psi_1), (\phi_2, \psi_2)$ are
distinct elements of $D$, then $\phi_1 \neq \phi_2$, and $\psi_1 \neq
\psi_2$.

\begin{lem}\label{diag}
Let $n \geq 0$, $\ep > 0$, and let $T : JH' \to JH$ be normalized.
Let
\begin{eqnarray*}
 A & = & \{(\phi,\psi) \in \levn\times\levn : \mbox{there exists\ }
                    \gamma,\xi \in \Gamma, \\
   & & \phi\in\gamma, \psi\in\xi, \mbox{such that\ } |\la
       TP'_n\gamma,P'_n\xi\ra| > \ep\}.
\end{eqnarray*}
Then $|D| \leq 1/\ep$ for all diagonal subsets $D$ of $A$.
\end{lem}

\begin{pf}
Let $D = \{(\phi_i,\psi_i): 1 \leq i \leq k\}$ be a diagonal subset of
$A$.  For each $i$, choose $\gamma_i, \xi_i \in \Gamma$ such that
$\phi_i \in \gamma_i, \psi_i \in \xi_i$, and $|\la
TP'_n\gamma_i,P'_n\xi_i\ra|
> \ep$.  Using the diagonality of $D$, we see that
$(P'_n\gamma_i)^k_{i=1}$ and $(P'_n\xi_i)^k_{i=1}$ are both
isometrically
equivalent to the $\ell^\infty(n)$-basis.
For $1 \leq i,j \leq k$, let $a_{ij} = \la
TP'_n\gamma_i,P'_n\xi_j\ra$.  Define $S, R: \ell^\infty(k) \to
\ell^1(k)$ by
$S(b_1,\ldots,b_k) = (\sum^k_{j=1}a_{ij}b_j)^k_{i=1}$ and
$R(b_1,\ldots,b_k) = (a_{ii}b_i)^k_{i=1}$ respectively.
Then $\|S\| \leq \|T\| = 1$, and $\|R\| \geq k\ep$.  However, by
\cite[Proposition 1.c.8]{LT}, $\|R\| \leq \|S\|$.  Therefore, $|D| = k
\leq 1/\ep$.
\end{pf}

\begin{lem}\label{longlem}
Let $(T_n)$ be a normalized weakly null sequence in $JH\epten JH$ such
that there is a sequence $0 = j_1 < j_2 < \cdots$
with $(P_{j_n} - P_{j_{n+1}})T_n(P_{j_n} - P_{j_{n+1}})' = T_n$ for
all
$n$. Then there is a subsequence $(T_{n_k})$ such that
\[ \sum|\la T_{n_k}\gamma,\xi\ra| < \infty \]
for all $\gamma, \xi \in \Gamma$.
\end{lem}

\begin{pf}
Note that $\|T_n\| = 1$ implies $|\la T_n\gamma,\xi\ra| \leq 1$ for
all $\gamma, \xi \in \Gamma$.  For all $m, n \in \N$, let
\begin{eqnarray*}
A(n,m) & = & \{(\phi,\psi)\in \levjn\times\levjn : \mbox{there exist
branches } \gamma \mbox{ and } \xi \mbox{ through }\\
& & \phi \mbox{ and } \psi \mbox{ respectively so that }
 |\la T_n\gamma,\xi\ra| > 1/2^m, \mbox{ and for all branches} \\
& & \gamma \mbox{ and $\xi$ through $\phi$ and
$\psi$ respectively, }
|\la T_n\gamma,\xi\ra| \leq 1/2^{m-1} \}.
\end{eqnarray*}
Fix $n$.  Let $B(n,1)$ be a maximal diagonal subset of $A(n,1)$.  Then
let
\[ C(n,1) = \bigcup_{(\phi,\psi)\in B(n,1)}\{(\phi_1,\psi_1): \phi_1 =
\phi \mbox{ or } \psi_1 = \psi\}. \]
Inductively, if $B(n,k)$ and $C(n,k)$ have been chosen for $k < m$,
let $B(n,m)$ be a maximal diagonal subset of $A(n,m)\backslash
\cup^{m-1}_{k=1}C(n,k)$. Then let
\[ C(n,m) = \bigcup_{(\phi,\psi)\in B(n,m)}\{(\phi_1,\psi_1): \phi_1 =
\phi \mbox{ or } \psi_1 = \psi\}. \]
It is easily seen that \\
(a) $\cup^\infty_{m=1}B(n,m)$ is a diagonal subset of
$\levjn\times\levjn$, \\
(b) for all $m$, $B(n,m) \subseteq A(n,m)\cap C(n,m)$, \\
(c) for all $k$, $\cup^k_{m=1}A(n,m) \subseteq \cup^k_{m=1}C(n,m)$.\\
In particular, by (b) and Lemma \ref{diag}, $|B(n,m)| \leq 2^m$ for
all $m, n$.  Also, if $(\phi,\psi) \in
\levjn\times\levjn\backslash\cup^k_{m=1}C(n,m)$, then (c) implies
$(\phi,\psi) \notin \cup^k_{m=1}A(n,m)$.  Hence if $\gamma, \xi \in
\Gamma$ pass through $\phi$ and $\psi$ respectively, then
\begin{equation}\label{size}
|\la T_n\gamma,\xi\ra| \leq 1/2^k.
\end{equation}
Now choose $N_1 \in \Pinf(\N)$ such that $|B(n,1)|$ is a constant, say
$b_1$, for all $n \in N_1$.  Write
\[ B(n,1) = \{(\phi(n,1,i),\psi(n,1,i)) : 1 \leq i \leq b_1\}  \]
for all $n \in N_1$. There exists $N'_1 \in \Pinf(N_1)$ such that for
each $i \leq b_1$,
$(\phi(n,1,i))_{n\in N'_1}$ as well as
$(\psi(n,1,i))_{n\in N'_1}$ are either strongly incomparable or
determine a branch.  Let
\[ I_1(\phi) =  \{1 \leq i \leq b_1: (\phi(n,1,i))_{n\in N'_1} \mbox{
determines a branch}. \} \]
Let the set of branches so determined be denoted by $\Gamma_1(\phi)$.
Define $I_1(\psi)$ and $\Gamma_1(\psi)$ similarly with regard to the
sequence $(\psi(n,1,i))_{n\in N'_1}$.  Since $(T_n)$ is weakly null,
so are $(T_n\gamma)$ and $(T'_n\xi)$ for all $\gamma, \xi \in \Gamma$.
Thus, because both $\Gamma_1(\phi)$ and $\Gamma_1(\psi)$ are finite
sets, Lemma \ref{seven} yields a set $N''_1 \in \Pinf(N'_1)$ such that
\[ \sup_{\xi\in\Gamma}\sum_{n\in N''_1}|\la T_n\gamma,\xi\ra| \leq 7
\spa \mbox{for all} \spa \gamma \in \Gamma_1(\phi) \]
and
\[ \sup_{\gamma\in\Gamma}\sum_{n\in N''_1}|\la\gamma,T'_n\xi\ra| \leq
7
\spa \mbox{for all} \spa \xi \in \Gamma_1(\psi). \]
Inductively, if $N''_m \in \Pinf(N)$ has been chosen, let $N_{m+1}
\in
\Pinf(N''_m)$ be such that $|B(n,m+1)|$ is a constant, say $b_{m+1}$,
for
all $n \in N_{m+1}$.  For $n \in N_{m+1}$, list
\[ B(n,m+1) = \{(\phi(n,m+1,i),\psi(n,m+1,i)) : 1 \leq i \leq
b_{m+1}\}, \]
choose $N'_{m+1} \in \Pinf(N_{m+1})$ such that for each $i \leq b_{m+1}$,
$(\phi(n,m+1,i))_{n\in N'_{m+1}}$ as well as
$(\psi(n,m+1,i))_{n\in N'_{m+1}}$ are either strongly incomparable or
determine a branch. Let
\[ I_{m+1}(\phi) =  \{1 \leq i \leq b_{m+1}: (\phi(n,m+1,i))_{n\in
N'_{m+1}} \mbox{ determines a branch}. \} \]
Let the set of branches so determined be denoted by
$\Gamma_{m+1}(\phi)$.
Define $I_{m+1}(\psi)$ and $\Gamma_{m+1}(\psi)$ similarly with regard
to the
sequence $(\psi(n,m+1,i))_{n\in N'_{m+1}}$.
If $\gamma \in \Gamma_{m+1}(\phi)$, $\xi \in \Gamma$, and $n \in
N'_{m+1}$, let $\phi$ and $\psi$ be the nodes of length $j_n$ in
$\gamma$ and $\xi$ respectively.  Then
$\phi = \phi(n,m+1,i)$ for some $i$.  But
$(\phi(n,m+1,i),\psi(n,m+1,i)) \in B(n,m+1)$.  Therefore, $(\phi,\psi)
\in C(n,m+1)$.  In particular, $(\phi,\psi) \notin \cup^m_{k=1}C(n,k)$.
By equation (\ref{size}), $|\la T_n\gamma,\xi\ra| \leq 1/2^m$.
Similarly, $|\la \gamma,T'_n\xi\ra| \leq 1/2^m$ for all $\gamma \in
\Gamma$ and $\xi \in \Gamma_{m+1}(\psi)$.  By
Lemma \ref{seven}, there is a set $N''_{m+1} \in \Pinf(N'_{m+1})$ such
that
\begin{equation}\label{supt}
\sup_{\xi\in\Gamma}\sum_{n\in N''_{m+1}}|\la T_n\gamma,\xi\ra| \leq
7/2^m
\spa \mbox{for all \spa $\gamma \in \Gamma_{m+1}(\phi)$},
\end{equation}
and
\begin{equation}\label{suptp}
\sup_{\gamma\in\Gamma}\sum_{n\in N''_{m+1}}|\la\gamma,T'_n\xi\ra|
\leq 7/2^m
\spa \mbox{for all \spa $\xi \in \Gamma_{m+1}(\psi)$}.
\end{equation}
Pick $n_1 < n_2 < \cdots$ such that $n_m \in N''_m$ for all $m$, and
let
\[ D(m) = \levjnm\times\levjnm\backslash\bigcup^m_{k=1}C(n_m,k). \]
For
all $m$, $\{C(n_m,1),\ldots,C(n_m,m),D(m)\}$ is a partition of
$\levjnm\times\levjnm$. Fix $\gamma, \xi \in \Gamma$.  We proceed to
estimate $\sum_m|\la T_{n_m}\gamma,\xi\ra|$.
For all $m \in \N$, let $\phi_m$ and $\psi_m$ be the nodes of length
$j_{n_m}$ in $\gamma$ and $\xi$ respectively.  Define $J_0 = \{m:
(\phi_m,\psi_m) \in D(m)\}$, and $J_k = \{m\geq k: (\phi_m,\psi_m) \in
C(n_m,k)\}$ for all $k \geq 1$.  Note that $\{J_0, J_1, J_2, \ldots\}$
is a partition of $\N$.  If $m \in J_0$, then equation (\ref{size})
yields $\tdot \leq 1/2^m$.  Consequently,
\begin{equation}\label{jnot}
\sum_{m\in J_0}\tdot \leq \sum_{m\in J_0}\frac{1}{2^m} \leq 1.
\end{equation}
Now fix $k \geq 1$ and $m \in J_k$. Then $m \geq k$ and hence
$n_m \in N_k$.  Thus $B(n_m,k)$ is listed as
$\{(\phi(n_m,k,i),\psi(n_m,k,i)): 1 \leq i \leq b_k\}$.  But
$(\phi_m,\psi_m) \in C(n_m,k)$.  Hence there exists $1 \leq i_m \leq
b_k$ such that either $\phi_m = \phi(n_m,k,i_m)$ or $\psi_m =
\psi(n_m,k,i_m)$.  Let $J_k(\phi) = \{m\in J_k: \phi_m =
\phi(n_m,k,i_m)\}$, and let $J_k(\psi) = J_k\backslash J_k(\phi)$.
Since $n_m \in N'_k$ as well, we may further subdivide these sets
into:
\begin{eqnarray*}
J_{k,1}(\phi) & = & \{m\in J_k(\phi): i_m \in I_k(\phi)\}, \\
J_{k,2}(\phi) & = & J_{k}(\phi)\backslash J_{k,1}(\phi), \\
J_{k,1}(\psi) & = & \{m\in J_k(\psi): i_m \in I_k(\psi)\}, \\
J_{k,2}(\psi) & = & J_{k}(\psi)\backslash J_{k,1}(\psi).
\end{eqnarray*}
Now $\{\phi(n_m,k,i_m): m \in J_{k,2}(\phi)\}$ is a subset of the
branch $\gamma$.  But it is also contained in the union of the
strongly incomparable sequences
$(\phi(n,k,i))_{n\in N'_{k}}$, $i \notin I_k(\phi)$.  Hence
$|J_{k,2}(\phi)| \leq 2b_k$.  Recalling equation (\ref{size}), we
obtain
\begin{equation}\label{jphitwo}
\sum_{m\in J_{k,2}(\phi)}\tdot \leq \sum_{m\in
J_{k,2}(\phi)}\frac{1}{2^{k-1}} \leq \frac{4b_k}{2^k}.
\end{equation}
Similarly,
\begin{equation}\label{jpsitwo}
\sum_{m\in J_{k,2}(\psi)}\tdot \leq \frac{4b_k}{2^k}.
\end{equation}
For any $m \in J_{k,1}(\phi)$, $\phi_m$ belongs to a branch in
$\Gamma_k(\phi)$. Let $\tilde{\gamma}$ be a branch in $\Gamma_k(\phi)$
such that $M(\tilde{\gamma}) = \{m\in
J_{k,1}(\phi):\phi_m\in\tilde{\gamma}\}$ is non-empty.  If
$M(\tilde{\gamma})$ is finite, let $m_0$ be its maximal element.  Then
$\phi_{m_0}$ belongs to both branches $\tilde{\gamma}$ and $\gamma$.
Therefore,
\begin{eqnarray*}
\sum_{m\in M(\tilde{\gamma})}\tdot & = &
\sum_{m\in M(\tilde{\gamma})\backslash\{m_0\}}\tdott + \tdotnot \\
& \leq & \sup_{\xi\in\Gamma}\sum_{n\in N''_k}\tdotnt + \tdotnot \\
& \leq & \frac{7}{2^{k-1}} + \frac{1}{2^{k-1}} = \frac{16}{2^k}
\end{eqnarray*}
by equations (\ref{supt}) and (\ref{size}) respectively.
On the other hand, if $M(\tilde{\gamma})$ is infinite, then
$\tilde{\gamma}$ and $\gamma$ coincide.  Thus the term containing
$T_{n_{m_0}}$ may simply be omitted from the above inequality.  Now
since $|\Gamma_k(\phi)| \leq b_k$, we obtain
\begin{equation}\label{jphione}
\sum_{m\in J_{k,1}(\phi)}\tdot \leq \frac{16b_k}{2^k}.
\end{equation}
Similarly,
\begin{equation}\label{jpsione}
\sum_{m\in J_{k,1}(\psi)}\tdot \leq \frac{16b_k}{2^k}.
\end{equation}
Combining equations (\ref{jnot})--(\ref{jpsione}), we see that
\[ \sum\tdot \leq \sum^\infty_{k=0}\sum_{m\in J_k}\tdot \leq 1 +
\sum^\infty_{k=1}\frac{40b_k}{2^k}. \]
To complete the proof, it remains to show that $\sum b_k/2^k <
\infty$.
Fix $m \in \N$.  By property (a) of the sets $B(n,k)$,
$\cup^m_{k=1}B(n_m,k)$ is a diagonal subset of
$\levjnm\times\levjnm$.  Hence $\{\phi(n_m,k,i): i\leq b_k,\ k\leq
m\}$
are all distinct, as are $\{\psi(n_m,k,i): i\leq b_k,\ k\leq m\}$.
By the definition of $B(n_m,k)$,  one can choose, for any $i\leq b_k,\
k\leq m$, branches $\gamma_{k,i}$ and $\xi_{k,i}$, passing through
$\phi(n_m,k,i)$ and $\psi(n_m,k,i)$ respectively, so that $\tdotki >
1/2^k$. Keeping in mind the assumption on $(T_n)$, we see that this
inequality remains valid if $\gamma_{k,i}$ and $\xi_{k,i}$ are
replaced by $\delta_{k,i} = P'_{j_{n_m}}\gamma_{k,i}$ and
$\zeta_{k,i} = P'_{j_{n_m}}\xi_{k,i}$
respectively. Since $(\delta_{k,i})$ and $(\zeta_{k,i})$ are
isometrically equivalent to the $\ell^\infty(\sum^m_{k=1}b_k)$-basis,
the map $S : \ell^\infty(\sum^m_{k=1}b_k) \to \ell^1(\sum^m_{k=1}b_k)$
\[ S(b_{k,i})^{b_k\ \,m}_{i=1 k=1} = (\sum_{k,i}\tdotdzp
b_{k,i})^{b_{k'}\ \,m}_{i'=1 k'=1} \]
has norm $\leq \|T_{n_m}\| = 1$.  Then \cite[Proposition 1.c.8]{LT}
implies that the ``diagonal'' of $S$ also has norm $\leq 1$.  But this
means
\[ \sum^m_{k=1}\frac{b_k}{2^k} \leq
\sum^m_{k=1}\sum^{b_k}_{i=1}|\tdotdz| \leq 1. \]
Since $m$ is arbitrary, $\sum b_k/2^k$ converges, as required.
\end{pf}

In order to apply Elton's extremal criterion, we need the following.
For convenience, we call an element of $JH'$ of the form $P'_m\gamma$,
where $m \geq 0$, and $\gamma \in \Gamma$,
a $m$-$\infty$ {\em segment}.

\begin{lem}\label{ipnset}
Let $W$ be the collection of all elements in $JH'$ of the form
$\sum^r_{i=1}a_iS_i$, where $r \in {\em \N}$, $\max|a_i| \leq 1$, and
there exist $m,\ n$, with $n$ possibly equal to $\infty$, so that
$\{S_1,\ldots,S_r\}$ is a set of pairwise disjoint $m$-$n$ segments.
Then $W$ is an i.p.n.\ subset of $JH'$.  Consequently, $W\times W$ is
an i.p.n.\ subset of $JH\epten JH$.
\end{lem}

\begin{pf}
The second assertion follows from the first by Lemma \ref{pn}.  If
$\sum^r_{i=1}a_iS_i \in W$, and $x \in JH$, then
\[ |\sum^r_{i=1}a_iS_ix| \leq \sum^r_{i=1}|S_ix| \leq \|x\|. \]
To complete the proof, it suffices to show that for all $x \in JH$,
there is a collection $\{S_1,\ldots,S_r\}$ of disjoint $m$-$n$
segments
($n$ possibly $= \infty$) such that
$\|x\| = \sum^r_{i=1}|S_ix|$.
Let $x \in JH$ be fixed.  For each $j$, choose an admissible
collection of $m_j$-$n_j$ segments $A_j$ such that
\begin{equation}\label{norming}
\|x\| = \lim_j\sum_{S\in A_j}|Sx|.
\end{equation}
If $(m_j)$ is unbounded, $\sum_{S\in
A_j}|Sx| \leq \|P_{m_j}x\| \to 0$ as $j \to \infty$.  Hence $x = 0$,
and the result is obvious.  If $(n_j)$ is bounded, then so is $(m_j)$.
Without loss of generality, we may assume that both $(m_j)$ and
$(n_j)$ are constant sequences with finite values, say, $m$ and $n$.
But as there are only finitely many sets of admissible $m$-$n$
segments,
the limiting value in (\ref{norming}) is attained, and the claim
holds. Finally, we consider the case when $(m_j)$ is bounded and
$(n_j)$ is unbounded.   Going to a subsequence, we may assume that
$(m_j)$ has a constant value, say $m$, and $n_j \to \infty$.  Then,
for each $j$, $|A_j| \leq 2^m$.  Using a subsequence again, we may
assume
that $|A_j| = r$ for some fixed $r$ for all $j$. For each $j$, write
$A_j = \{S_1(j),\ldots,S_r(j)\}.$  Choose a subsequence $(j_k)$ such
that $(S_i(j_k))_k$ converges weak* to some $S_i$ for every $1 \leq i
\leq r$.  It is easy to see that $\{S_1,\ldots,S_r\}$ is a collection
of pairwise disjoint $m$-$\infty$ segments.  From equation
(\ref{norming}), we deduce that $\|x\| = \sum^r_{i=1}|S_ix|$, as
desired.
\end{pf}

\begin{lem}\label{mainlem}
Let $(T_n)$ be as in Lemma \ref{longlem}, then $[T_n]$ contains a copy
of $\cz$.
\end{lem}

\begin{pf}
Let $W$ be as in Lemma \ref{ipnset}.  Choose a subsequence $(T_{n_k})$
as given by Lemma \ref{longlem}.  Then we have $\sum|\la
T_{n_k}w,v\ra| < \infty$ for all $(w,v) \in W\times W$.  But by Lemma
\ref{ipnset}, $W\times W$ is an i.p.n.\ subset of $JH\epten JH$.
Hence
Elton's extremal criterion \cite{E} assures us that $[T_{n_k}]$
contains a copy of $\cz$.
\end{pf}

\begin{thm}
The space $JH\epten JH$ is $\cz$-saturated.
\end{thm}

\begin{pf}
For all $m \in \N$, let
\[ E_m = \kw(JH',(1-P_m)JH) \spa \mbox{and} \spa F_m =
\kw((1-P_m)'JH',JH). \]
Both $E_m$ and $F_m$ are isomorphic to a direct
sum of a finite number of copies of $JH$, hence they are
$\cz$-saturated by Lemma \ref{EF}.
Let $(S_n)$ be a normalized basic sequence in $JH\epten JH$.  If there
exist
a subsequence $(R_n)$ of $(S_n)$ and $m \in \N$\ such that $(R_n)$
is dominated by $((1-P_m)R_n\oplus R_n(1-P_m)') \in E_m\oplus F_m$,
then $(R_n)$ is equivalent to a sequence in $E_m\oplus F_m$, which is
a $\cz$-saturated space by Lemma \ref{EF}.  Thus $[R_n]$, and
consequently $[S_n]$, contains a copy of $\cz$.  Otherwise, for all $m
\in \N$\ and every subsequence $(R_n)$ of $S_n$,
\[ \inf\left\{\|\sum_na_n(1-P_m)R_n\| + \|\sum_na_nR_n(1-P_m)'\| :
(a_n)
\in \czz, \|\sum a_nR_n\| = 1\right\} = 0. \]
Also, it is clear that for any $T \in JH\epten JH$,
$\lim_n(1-P_n)T(1-P_n)' =
T$ in norm. Using these observations and a
standard perturbation argument, we obtain a normalized block basis
of $(S_n)$
which is
equivalent to some sequence $(T_n)$ satisfying the hypotheses of Lemma
\ref{longlem}. By Lemma \ref{mainlem}, $[T_n]$ contains a copy of
$\cz$.  Thus, so does $[S_n]$.
\end{pf}

Hagler \cite{H} proved that, in fact, $JH$ has {\em property}\/ (S):
every normalized weakly null sequence has a $\cz$-subsequence.
Coupled with the absence of $\ell^1$, property (S) implies
$\cz$-saturation. In \cite{KO}, it was also shown that property (S)
implies property (u).  Thus, we may ask\\

\noindent{\em Question}:  Does $JH\epten JH$ has property (S) or
property (u)?

\flushleft
\vspace{.2in}
Department of Mathematics\\National University of Singapore\\
Singapore 0511\\ E-mail(bitnet) : matlhh@nusvm


\baselineskip 3ex

\begin{thebibliography}{99}

\bibitem{E} J. Elton, {\em Extremely weakly unconditionally
convergent series}, Israel J.\ Math.\ {\bf 40}(1981), 255-258.

\bibitem{F1} V. F. Fonf, {\em On a property of Lindenstrauss-Phelps
spaces}, Funct.\ Anal.\ Appl.\ {\bf 13}(1979), 79-80 (translated from
Russian).

\bibitem{F2} V. F. Fonf, {\em Polyhedral Banach spaces},
Matematicheskie
Zametki {\bf 30}(1981), 627-634 (translated from Russian).

\bibitem{H} James Hagler, {\em A counterexample to several questions
about Banach spaces }, Studia Math.\ {\bf 60}(1977), 289-308.

\bibitem{KO} H. Knaust and E. Odell, {\em On $\cz$-sequences in Banach
spaces}, Israel J.\ Math.\ {\bf67}(1989), 153-169.

\bibitem{L} Denny H. Leung, {\em Embedding $\ell^1$ into tensor
products of Banach spaces}, Functional Analysis (eds.\ E. Odell and H.
Rosenthal), Lecture Notes in Math., vol.\ 1470, Springer-Verlag,
Berlin,
1991, 171-176.

\bibitem{LT} Joram Lindenstrauss and Lior Tzafriri, ``Classical Banach
Spaces I, Sequence Spaces'', Springer-Verlag, Berlin, 1977.

\bibitem{Ro} H. Rosenthal, {\em A characterization of Banach spaces
containing $\ell^1$}, Proc.\ Nat.\ Acad.\ Sci.\ (U.S.A.), {\bf
71}(1974),
2411-2413.

\bibitem{R1} H. Rosenthal, Class Notes, Topics course in analysis,
University of Texas at Austin.

\bibitem{R} H. Rosenthal, {\em Some aspects of the subspace structure
of infinite dimensional Banach spaces}, Approximation Theory and
Functional Analysis (ed.\ C. Chuy), Academic Press, 1990.

\bibitem{Ru} W. Ruess, {\em Duality and geometry of spaces of compact
operators}, Functional Analysis: Surveys and Recent Results III, (eds.\
K.-D. Bierstedt and B. Fuchssteiner) Math.\ Studies, no.\ 90,
North-Holland, Amsterdam, 1984, 59-78.

\bibitem{S} H.\ H.\ Schaefer, ``Banach Lattices and Positive
Operators'', Springer-Verlag, Berlin, 1974.

\bibitem{Sem} Z.\ Semadeni, ``Banach Spaces of Continuous Functions
I'', Monografie Matematyczne, PWN-Polish Scientific Publishers 55,
Warszawa, 1971.

\end{thebibliography}
\end{document}